\documentclass[12pt,a4paper,oneside]{article}
\usepackage{amsmath,amsfonts,amssymb,latexsym}

\newtheorem{theorem}{Theorem}[section]
\newtheorem{lemma}[theorem]{Lemma}

\newenvironment{proof}
{\bigskip\noindent{\sc Proof.}\ \ \rm }{\hfill$\Box$\bigskip}
\title{ On a new version of the It\^o's formula\\ for the stochastic heat equation }
\author{ Alberto Lanconelli }
\date{\empty}
\begin{document}
\maketitle { \noindent \small
\begin{tabular}{ll}
& Dipartimento di Matematica\\
& Universita' degli Studi di Bari\\
& Via E. Orabona, 4\\
& 70125 Bari - Italia\\

& Email: lanconelli@dm.uniba.it \\

\end{tabular}
} \numberwithin{equation}{section}

\bigskip

\begin{abstract}
We derive an It\^o's-type formula for the one dimensional stochastic
heat equation driven by a space-time white noise. The proof is based
on elementary properties of the $\mathcal{S}$-transform and on the
explicit representation of the solution process. We also discuss the
relationship with other versions of this It\^o's-type formula
existing in literature.
\end{abstract}
\section{Introduction}
Consider the following stochastic partial differential equation
(SPDE):
\begin{eqnarray}
\Big{ \{ }
\begin{array}{ll}
\partial_tu_t(x)=\frac{1}{2}\partial_{xx} u_t(x)+\dot{W}_{t,x} \\
u_t(0)=u_t(1)=0,\quad u_0=0
\end{array}
\end{eqnarray}
where $\partial_t:=\frac{\partial}{\partial t}$,
$\partial_{xx}:=\frac{\partial^2}{\partial x^2}$,
$\dot{W}_{t,x}:=\frac{\partial^2W_{t,x}}{\partial t\partial x}$ and
$\{W_{t,x}, t\in [0,T], x\in [0,1]\}$ is a Brownian sheet. By
solution to this equation we mean an adapted two parameter
stochastic process $\{u_t(x), t\in [0,T], x\in [0,1]\}$ such that
for $t\in [0,T]$ and $x\in [0,1]$,
$$
u_t(0)=u_t(1)=0,\quad u_0(x)=0,
$$
and such that for all $l\in C_0^2(]0,1[)$ the following equality
$$
\langle u_t,l\rangle=\langle u_0,l\rangle+\frac{1}{2}\int_0^t\langle
u_s,l''\rangle ds+\int_0^t\langle l,dW_s\rangle,
$$
holds for $t\in [0,T]$. Here $\langle,\rangle$ denotes the inner
product in $\mathcal{L}^2([0,1])$ and
$$
\int_0^t\langle l,dW_s\rangle:=\int_0^t\int_0^1l(x)dW_{s,x}.
$$
It is well known (see e.g. \cite{W}) that equation (1.1) has a
unique solution $\{u_t(x), t\in [0,T], x\in [0,1]\}$ which is
continuous in the variables $(t,x)$ and that it can be represented
as
$$
u_t(x)=\int_0^t\int_0^1g_{t-s}(x,y)dW_{s,y}
$$
where $\{g_t(x,y), t\in [0,T], x,y\in [0,1]\}$ is the fundamental
solution of the heat equation with homogenous Dirichlet boundary
conditions, i.e.
\begin{eqnarray}
\Big{ \{ }
\begin{array}{ll}
\partial_tg_t(x,y)=\frac{1}{2}\partial_{xx} g_t(x,y)\\
g_t(0,y)=g_t(1,y)=0,\quad g_0(x,y)=\delta(x-y)
\end{array}
\end{eqnarray}
Since for fixed $x\in [0,1]$ the process $t\mapsto u_t(x)$ is not a
semimartingale, we can not deal with it by means of the classical
stochastic calculus. It is therefore natural to ask whether an
It\^o's-type formula can be found for this kind of process. Two
recent papers, \cite{GNT} and \cite{Z}, are devoted to the
investigation of this problem. In \cite{GNT} the authors develop a
Malliavin calculus for the solution process $u_t(x)$ in order to
obtain an It\^o's-type formula whose proof makes also use of
projections on Wiener chaoses of different orders. In \cite{Z} the
author applies the ordinary It\^o's formula to a regularized version
of the solution process $u_t(x)$; then he studies the limit when
that regularized process converges to the real one. The main feature
of this procedure is the appearance of a renormalization of the
square of an infinite
dimensional stochastic distribution.\\
The aim of the present paper is to propose an alternative approach
(and a corresponding version of the It\^o's-type formula) to the
above mentioned problem which is somehow in between the techniques
of the articles \cite{GNT} and \cite{Z}; in fact we utilize notions
of white noise analysis, an infinite dimensional stochastic
distributions theory, and prove the resulting formula via scalar
products with test functions, more precisely stochastic
exponentials. The key point is the gaussianity of the solution
process and its explicit representation as a stochastic convolution.
Our formula is in the spirit of the It\^o's-type formula for
gaussian processes derived in the paper \cite{NT}. Moreover the idea
of using the properties of the semigroup associated to the one
dimensional Brownian motion (see the proof of Theorem 2.2) is
analogous to the one used
 in \cite{L3}.\\
The paper is structured as follows: in Section 2.1 we recall basic
notions and facts from the white noise theory; then in Section 2.2
we state and prove our main result; finally in the concluding
section we compare our formula with those already existing in the
literature.
\section{Main result}
\subsection{Preliminaries}
In this section we fix the notation and recall some basic results
from the white noise theory. For additional information about this
topic we refer the reader to the books \cite{HOUZ} and \cite{Kuo}
and to the paper \cite{DPV}.\\
Let $(\Omega,\mathcal{F},\mathcal{P})$ be a complete
probability space and $\{W_{t,x}, t\in[0,T], x\in [0,1]\}$ a
Brownian sheet defined on it. Assume that
$$
\mathcal{F}=\sigma(W_{t,x}, t\in[0,T], x\in [0,1]),
$$
so that each element
$X\in\mathcal{L}^2(\Omega,\mathcal{F},\mathcal{P})$ can be
decomposed as a sum of multiple It\^o integrals w.r.t. the Brownian
sheet $W$, i.e.
$$
X=\sum_{n\geq 0}I_n(h_n)\mbox{ (convergence in
}\mathcal{L}^2(\Omega,\mathcal{F},\mathcal{P})),
$$
where for $n\geq 0$, $h_n$ is a deterministic function belonging to
$\mathcal{L}^2(([0,T]\times [0,1])^n)$ and $I_n(h_n)$ is the $n$-th
order multiple It\^o integral of $h_n$.\\
For example, the solution of the SPDE (1.1) has the decomposition:
$$
u_t(x)=\int_0^t\int_0^1g_{t-s}(x,y)dW_{s,y}=I_1(1_{[0,t]}(\cdot)g_{t-\cdot}(x,\cdot)).
$$
If $A:D(A)\subset\mathcal{L}^2([0,T]\times
[0,1])\to\mathcal{L}^2([0,T]\times [0,1])$ is the unbounded operator
$$
Ah(t,x):=A_tA_xh(t,x):=\sqrt{-\frac{\partial^2}{\partial
t^2}}\sqrt{-\frac{\partial^2}{\partial x^2}}h(t,x)
$$
with periodic boundary conditions, we define its \emph{second
quantization} as the operator,
\begin{eqnarray*}
\Gamma(A):D(\Gamma(A))\subset\mathcal{L}^2(\Omega,\mathcal{F},\mathcal{P})&\to&
\mathcal{L}^2(\Omega,\mathcal{F},\mathcal{P})\\
\sum_{n\geq 0}I_n(h_n)&\mapsto&\Gamma(A)(\sum_{n\geq
0}I_n(h_n)):=\sum_{n\geq 0}I_n(A^{\otimes n}h_n).
\end{eqnarray*}
For $p\geq 1$ the space
$$
(S)_p:=\{X\in\mathcal{L}^2(\Omega,\mathcal{F},\mathcal{P})\mbox{
s.t. }E[|\Gamma(A^p)X|^2]<+\infty\}
$$
is a subset of $\mathcal{L}^2(\Omega,\mathcal{F},\mathcal{P})$ and
if $q>p$ then $(S)_q\subset (S)_p$. The \emph{Hida's test functions
space} is defined as
$$
(S):=\bigcap_{p\geq 1}(S)_p.
$$
Its dual w.r.t. the inner product of
$\mathcal{L}^2(\Omega,\mathcal{F},\mathcal{P})$ is called
\emph{Hida's distributions space} and denoted by $(S)^*$. It can be
shown that
$$
(S)^*=\bigcup_{p\geq 1}(S)_{-p};
$$
moreover by construction
$$
(S)\subset\mathcal{L}^2(\Omega,\mathcal{F},\mathcal{P})\subset
(S)^*.
$$
For $f\in C_0^{\infty}(]0,T[\times ]0,1[)$, the random variable
$$
\mathcal{E}_T(f):=\exp\{\int_0^T\int_0^1f(s,y)dW_{s,y}-\frac{1}{2}\int_0^T\int_0^1f^2(s,y)dyds\},
$$
belongs to $(S)$; therefore if $\langle\langle,\rangle\rangle$
denotes the dual pairing between $(S)^*$ and $(S)$ then the
application
$$
X\in (S)^*\mapsto\mathcal{S}(X)(f):=\langle\langle
X,\mathcal{E}_T(f)\rangle\rangle
$$
is well defined and it is called $\mathcal{S}$\emph{-transform}. The
quantity $\mathcal{S}(X)(f), f\in C_0^{\infty}(]0,T[\times ]0,1[)$
identifies uniquely the Hida's distribution $X$. In particular,
given $X,Y\in (S)^*$, we denote by $X\diamond Y$ the unique element
of $(S)^*$ such that
$$
\mathcal{S}(X\diamond
Y)(f)=\mathcal{S}(X)(f)\mathcal{S}(Y)(f),\mbox{ for all }f\in
C_0^{\infty}(]0,T[\times ]0,1[);
$$
the stochastic distribution $X\diamond Y$ is named \emph{Wick
product} of $X$ and $Y$.\\
Notice also that if
$X\in\mathcal{L}^2(\Omega,\mathcal{F},\mathcal{P})$ then
$$
\mathcal{S}(X)(f)=\langle\langle
X,\mathcal{E}_T(f)\rangle\rangle=E[X\mathcal{E}_T(f)],
$$
and if $\xi$ is an It\^o integrable stochastic process then
$$
\mathcal{S}(\int_0^T\int_0^1\xi_{t,x}dW_{t,x})(f)=E[\int_0^T\int_0^1\xi_{t,x}dW_{t,x}\mathcal{E}_T(f)]=
\int_0^T\int_0^1E[\xi_{t,x}\mathcal{E}_T(f)]dxdt.
$$
\subsection{It\^o's-type formula}
Before the main theorem of this paper, we state and prove the
following auxiliary result.
\begin{lemma}
Let $u_t(x)$ be the solution of the SPDE (1.1); then
$\partial_{xx}u_t(x)\in (S)^*$.
\end{lemma}
\begin{proof}
It is known that the fundamental solution $g_t(x,y)$ of the heat
equation with homogenous Dirichlet boundary conditions (1.2) can be
represented as
$$
g_t(x,y)=\sum_{n\geq 1}e^{-\lambda_nt}e_n(x)e_n(y),
$$
where $e_n(x)=\sqrt{2}sin(\pi nx)$ and $\lambda_n=\pi^2n^2$. This
gives
\begin{eqnarray*}
\partial_{xx}g_t(x,y)&=&\sum_{n\geq
1}e^{-\lambda_nt}(-\lambda_n)e_n(x)e_n(y)\\
&=&\sum_{n\geq 1}e^{-\lambda_nt}e_n(x)\partial_{yy}e_n(y)\\
&=&\sum_{n\geq 1}e^{-\lambda_nt}e_n(x)(-A^2_ye_n(y)).
\end{eqnarray*}
Therefore,
\begin{eqnarray*}
E[|\Gamma(A^{-2})\partial_{xx}u_t(x)|^2]&=&\int_0^t\int_0^1|A_s^{-2}A_y^{-2}\partial_{xx}g_{t-s}(x,y)|^2dyds\\
&=&\int_0^t\int_0^1|A_s^{-2}g_{t-s}(x,y)|^2dyds\\
&\leq&\int_0^t\int_0^1|g_{t-s}(x,y)|^2dyds<+\infty,
\end{eqnarray*}
proving that $\partial_{xx}u_t(x)\in (S)_{-2}\subset (S)^*$.
\end{proof}
\begin{theorem}
For any $\varphi\in C^2_b(\mathbb{R})$ and $l\in C_0^2(]0,1[)$ one
has
\begin{eqnarray}
\langle\varphi(u_t),l\rangle&=&\langle\varphi(u_0),l\rangle+\int_0^t\langle\varphi'(u_s)l,dW_{s}\rangle
+\frac{1}{2}\langle\int_0^t\varphi'(u_s)\diamond\partial_{xx}u_sds,l\rangle\nonumber\\
&&+\frac{1}{2}\langle\int_0^t\varphi''(u_s)d\sigma^2(s;\cdot),l\rangle,
\end{eqnarray}
where
\begin{itemize}
\item $\sigma^2(s,x)=E[|u_s(x)|^2]=\int_0^s\int_0^1g_{s-v}^2(x,y)dydv$;
\item $\int_0^t\varphi''(u_s(x))d\sigma^2(s;x)$ is a Stiltjes integral in
the variable $s$;
\item $\varphi'(u_s)\diamond\partial_{xx}u_s$ is the Wick
product between $\varphi'(u_s)$ and $\partial_{xx}u_s$.
\end{itemize}
\end{theorem}
\begin{proof}
We aim at proving that
$$
E[\langle\varphi(u_t),l\rangle\mathcal{E}_T(f)]=E[\mathcal{R}\mathcal{E}_T(f)],
$$
for all $f\in C_0^{\infty}(]0,T[\times ]0,1[)$ where $\mathcal{R}$
denotes the r.h.s. of equation (2.1). This fact together with the
property
$$
span\{\mathcal{E}_T(f), f\in C_0^{\infty}(]0,T[\times ]0,1[)\}\mbox{
is dense in }\mathcal{L}^2(\Omega,\mathcal{F},\mathcal{P})
$$
will imply that
$$
\langle\varphi(u_t),l\rangle=\mathcal{R}\quad\mathcal{P}-\mbox{almost
surely},
$$
which is the statement of the theorem.\\
Let us fix an arbitrary
$f\in C_0^{\infty}(]0,T[\times ]0,1[)$; a simple application of the
Girsanov's theorem yields:
\begin{eqnarray*}
E[\varphi(u_t(x))\mathcal{E}_T(f)]&=&E\Big[\varphi\Big(\int_0^t\int_0^1g_{t-s}(x,y)dW_{s,y}\Big)\mathcal{E}_T(f)\Big]\\
&=&E\Big[\varphi\Big(\int_0^t\int_0^1g_{t-s}(x,y)dW_{s,y}+\int_0^t\int_0^1g_{t-s}(x,y)f(s,y)dyds\Big)\Big].\\
\end{eqnarray*}
Now observe that
$$
\int_0^t\int_0^1g_{t-s}(x,y)dW_{s,y}+\int_0^t\int_0^1g_{t-s}(x,y)f(s,y)dyds,
$$
is a Gaussian random variable with mean given by
$$
m(t,x):=\int_0^t\int_0^1g_{t-s}(x,y)f(s,y)dyds,
$$
and variance equal to
$$
\sigma^2(t,x):=\int_0^t\int_0^1g^2_{t-s}(x,y)dyds.
$$
Therefore we can write
\begin{eqnarray}
E[\varphi(u_t(x))\mathcal{E}_T(f)]&=&E\Big[\varphi\Big(
\int_0^t\!\!\int_0^1g_{t-s}(x,y)dW_{s,y}+\int_0^t\!\!\int_0^1g_{t-s}(x,y)f(s,y)dyds\Big)\Big]\nonumber\\
&=&(P_{\sigma^2(t,x)}\varphi)(m(t,x)),
\end{eqnarray}
where $\{P_t\}_{t\geq 0}$ denotes the one dimensional heat
semigroup,
$$
(P_t\varphi)(x):=\int_{\mathbb{R}}\varphi(y)\frac{1}{\sqrt{2\pi
t}}e^{-\frac{(x-y)^2}{2t}}dy.
$$
Identity (2.2) turn out to be very useful; in fact we can apply the
ordinary chain rule and get
\begin{eqnarray*}
E[\varphi(u_t(x))\mathcal{E}_T(f)]&=&(P_{\sigma^2(t,x)}\varphi)(m(t,x))\\
&=&\varphi(u_0(x))+\int_0^t\frac{d}{dv}(P_{\sigma^2(v,x)}\varphi)(m(v,x))dv\\
&=&\varphi(u_0(x))+\int_0^t(P_{\sigma^2(v,x)}\varphi')(m(v,x))\frac{dm(v,x)}{dv}dv\\
&+&\int_0^t\frac{1}{2}(P_{\sigma^2(v,x)}\varphi'')(m(v,x))d\sigma^2(v;x)\\
&=:&\varphi(u_0(x))+\mathcal{I}(t,x)+\mathcal{II}(t,x).
\end{eqnarray*}
Recalling that
$$
m(v,x)=\int_0^v\int_0^1g_{v-s}(x,y)f(s,y)dyds,
$$
we have
\begin{eqnarray*}
\frac{dm(v,x)}{dv}&=&f(v,x)+\int_0^v\int_0^1\frac{1}{2}\partial_{xx}g_{v-s}(x,y)f(s,y)dyds\\
&=&f(v,x)+\frac{1}{2}\partial_{xx}\Big(\int_0^v\int_0^1g_{v-s}(x,y)f(s,y)dyds\Big),
\end{eqnarray*}
and hence
\begin{eqnarray*}
\mathcal{I}(t,x)&=&\int_0^t(P_{\sigma^2(v,x)}\varphi')(m(v,x))\Big(f(v,x)+
\frac{1}{2}\partial_{xx}\Big(\int_0^v\int_0^1g_{v-s}(x,y)f(s,y)dyds\Big)\Big)dv\\
&=&\int_0^t(P_{\sigma^2(v,x)}\varphi')(m(v,x))f(v,x)dv\\
&+&\int_0^t(P_{\sigma^2(v,x)}\varphi')(m(v,x))
\frac{1}{2}\partial_{xx}\Big(\int_0^v\int_0^1g_{v-s}(x,y)f(s,y)dyds\Big)dv\\
&=&\int_0^tE[\varphi'(u_v(x))\mathcal{E}_T(f)]f(v,x)dv+\int_0^tE[\varphi'(u_v(x))\mathcal{E}_T(f)]
\frac{1}{2}\partial_{xx}E[u_v(x)\mathcal{E}_T(f)]dv,
\end{eqnarray*}
where the last equality is due to identities (2.2). Moreover since
from Lemma 2.1 $\partial_{xx}u_v(x)\in (S)^*$ we have
$$
\partial_{xx}E[u_v(x)\mathcal{E}_T(f)]=\partial_{xx}\mathcal{S}(u_v(x))(f)=\mathcal{S}(\partial_{xx}u_v(x))(f),
$$
where $\mathcal{S}(\partial_{xx}u_v(x))(f)$ denotes the
$\mathcal{S}$-transform of $\partial_{xx}u_v(x)$. \\
If now $l\in C_0^2(]0,1[)$ is a test function in the space variable
$x$, by the properties of the Wick product we conclude that
\begin{eqnarray*}
\int_0^1\mathcal{I}(t,x)l(x)dx&=&\int_0^1\int_0^tE[\varphi'(u_v(x))\mathcal{E}_T(f)]l(x)f(v,x)dvdx\\
&+&\int_0^1\int_0^tE[\varphi'(u_v)\mathcal{E}_T(f)]
\frac{1}{2}\mathcal{S}(\partial_{xx}u_v(x))(f)l(x)dvdx\\
&=&\int_0^1\int_0^tE[\varphi'(u_v(x))\mathcal{E}_T(f)]l(x)f(v,x)dvdx\\
&+&\int_0^1\int_0^t\mathcal{S}(\varphi'(u_v(x))\diamond\frac{1}{2}\partial_{xx}u_v(x))(f)l(x)dvdx\\
&=&E\Big[\Big(\int_0^t\int_0^1\varphi'(u_v(x))l(x)dW_{v,x}\\
&+&\int_0^1\Big(\int_0^t\varphi'(u_v(x))\diamond\frac{1}{2}\partial_{xx}u_v(x)dv\Big)l(x)dx\Big)\mathcal{E}_T(f)\Big].
\end{eqnarray*}
Looking again at identities (2.2) we also discover that
$$
\mathcal{II}(t,x)=\int_0^t\frac{1}{2}E[\varphi''(u_s(x))\mathcal{E}_T(f)]d\sigma^2(s;x);
$$
combining the expressions for $\mathcal{I}$ and $\mathcal{II}$ we
can now conclude that
\begin{eqnarray*}
\langle
E[\varphi(u_t)\mathcal{E}_T(f)],l\rangle&=&\langle\varphi(u_0),l\rangle\\
&+&E\Big[\Big(\int_0^t\int_0^1\varphi'(u_v(x))l(x)dW_{v,x}+
\langle\int_0^t\varphi'(u_v)\diamond\frac{1}{2}\partial_{xx}u_vdv,l\rangle\Big)\mathcal{E}_T(f)\Big]\\
&+&\langle\int_0^t\frac{1}{2}E[\varphi''(u_s)\mathcal{E}_T(f)]d\sigma^2(v;\cdot),l\rangle\\
&=&\langle\varphi(u_0),l\rangle\\
&+&E\Big[\Big(\int_0^t\langle\varphi'(u_v)l,dW_{v}\rangle+
\langle\int_0^t\varphi'(u_v)\diamond\frac{1}{2}\partial_{xx}u_vdv,l\rangle\Big)\mathcal{E}_T(f)\Big]\\
&+&E\Big[\Big(\langle\frac{1}{2}\int_0^t\varphi''(u_s)d\sigma^2(v;\cdot),l\rangle\Big)\mathcal{E}_T(f)\Big].
\end{eqnarray*}
This completes the proof.
\end{proof}
\subsection{Comparisons}
\begin{description}
\item[Zambotti's formula]:
In \cite{Z} the author obtains the following It\^o's formula:
\begin{eqnarray}
\langle\varphi(u_t),l\rangle&=&\langle\varphi(u_0),l\rangle+\frac{1}{2}\int_0^t\langle
l'',\varphi(u_s)\rangle
ds+\int_0^t\langle\varphi'(u_s)l,dW_s\rangle\nonumber\\
&&-\frac{1}{2}\int_0^t\langle l,:\big|\partial_x
u_s\big|^2:\varphi''(u_s)\rangle ds,
\end{eqnarray}
where the last term is defined as the limit of renormalized
diverging quantities. The procedure to derive this formula is to
approximate the solution of the SPDE via the smoother process
$$
u_t^{\epsilon}(x):=\int_0^t\int_0^1g_{t-s+\epsilon}(x,y)dW_{s,y},
$$
and then to pass to the limit for $\epsilon\to 0$. We are now going
to show how to manipulate formula (2.1) to make it look like
(2.3).\\
Following the same line of reasoning explained in the proof of
Theorem 2.2 with $u_t^{\epsilon}(x)$ instead of $u_t(x)$, one can
easily see that everything can be carried in a similar manner; the
only difference consists in the fact that in this case the function
$$
\sigma_{\epsilon}^2(t,x):=E[|u_t^{\epsilon}(x)|^2]=\int_0^t\int_0^1g^2_{t-s+\epsilon}(x,y)dyds,
$$
is differentiable w.r.t. $t$; therefore the last term of formula
(2.1) can be rewritten as
\begin{eqnarray*}
\frac{1}{2}\int_0^t\varphi''(u_s^{\epsilon}(x))d\sigma_{\epsilon}^2(s;x)&=&
\frac{1}{2}\int_0^t\varphi''(u_s^{\epsilon}(x))\frac{d\sigma_{\epsilon}^2(s,x)}{ds}ds\\
=\frac{1}{2}\int_0^t\varphi''(u_s^{\epsilon}(x))\Big(\int_0^1g_{\epsilon}^2(x,y)dy&+&
\int_0^s\int_0^1\partial_s(g^2_{s-v+\epsilon}(x,y))dydv\Big)ds\\
=\frac{1}{2}\int_0^t\varphi''(u_s^{\epsilon}(x))\Big(\int_0^1g_{\epsilon}^2(x,y)dy&+&
\int_0^s\int_0^1g_{s-v+\epsilon}(x,y)\partial_{xx}g_{s-v+\epsilon}(x,y)dydv\Big)ds\\
=\frac{1}{2}\int_0^t\varphi''(u_s^{\epsilon}(x))\int_0^1\!\!g_{\epsilon}^2(x,y)dyds&+&
\frac{1}{2}\int_0^t\!\!\varphi''(u_s^{\epsilon}(x))
\int_0^s\!\!\int_0^1\!\!g_{s-v+\epsilon}(x,y)\partial_{xx}g_{s-v+\epsilon}(x,y)dydvds\\
=\frac{1}{2}\int_0^t\varphi''(u_s^{\epsilon}(x))\int_0^1g_{\epsilon}^2(x,y)dyds&+&
\frac{1}{2}\int_0^t\int_0^s\int_0^1D_{v,y}\varphi'(u^{\epsilon}_s(x))D_{v,y}\partial_{xx}u_s^{\epsilon}(x)dydvds\\
=\frac{1}{2}\int_0^t\varphi''(u_s^{\epsilon}(x))\int_0^1g_{\epsilon}^2(x,y)dyds
&+&\frac{1}{2}\int_0^t\varphi'(u^{\epsilon}_s(x))\partial_{xx}u_s^{\epsilon}(x)-
\varphi'(u^{\epsilon}_s(x))\diamond\partial_{xx}u_s^{\epsilon}(x)ds.
\end{eqnarray*}
Here we used a property of the Wick product which shows its
interplay with the Hida-Malliavin derivative, namely
$$
X\diamond
\int_0^T\int_0^1h(s,y)dW_{s,y}=X\int_0^T\int_0^1h(s,y)dW_{s,y}-\int_0^T\int_0^1D_{s,y}Xh(s,y)dyds.
$$
See \cite{Kuo} and \cite{N} for details. Therefore equation (2.1)
becomes
\begin{eqnarray*}
\langle\varphi(u^{\epsilon}_t),l\rangle&=&\langle\varphi(u^{\epsilon}_0),l\rangle+\int_0^t\langle\varphi'(u^{\epsilon}_s)l,dW_{s}\rangle
+\frac{1}{2}\langle\int_0^t\varphi'(u^{\epsilon}_s)\diamond\partial_{xx}u^{\epsilon}_sds,l\rangle\\
&&+\frac{1}{2}\langle\int_0^t\varphi''(u_s^{\epsilon})\int_0^1g_{\epsilon}^2(\cdot,y)dyds,l\rangle\\
&&+\frac{1}{2}\langle\int_0^t\varphi'(u^{\epsilon}_s)\partial_{xx}u_s^{\epsilon}-
\varphi'(u^{\epsilon}_s)\diamond\partial_{xx}u_s^{\epsilon}ds,l\rangle\\
&=&\langle\varphi(u^{\epsilon}_0),l\rangle+\int_0^t\langle\varphi'(u^{\epsilon}_s)l,dW_{s}\rangle\\
&&+\frac{1}{2}\langle\int_0^t\varphi''(u_s^{\epsilon})\int_0^1g_{\epsilon}^2(\cdot,y)dyds,l\rangle
+\frac{1}{2}\langle\int_0^t\varphi'(u^{\epsilon}_s)\partial_{xx}u_s^{\epsilon}ds,l\rangle.
\end{eqnarray*}
Moreover
$$
\varphi'(u^{\epsilon}_s(x))\partial_{xx}u_s^{\epsilon}(x)=\partial_{xx}\varphi(u_s^{\epsilon}(x))-
\varphi''(u^{\epsilon}_s(x))(\partial_{x}u_s^{\epsilon}(x))^2;
$$
a substitution in the previous equation gives
\begin{eqnarray*}
\langle\varphi(u^{\epsilon}_t),l\rangle&=&
\langle\varphi(u^{\epsilon}_0),l\rangle+\int_0^t\langle\varphi'(u^{\epsilon}_s)l,dW_{s}\rangle
+\frac{1}{2}\langle\int_0^t\varphi''(u_s^{\epsilon})(\int_0^1g_{\epsilon}^2(\cdot,y)dy)ds,l\rangle\\
&&+\frac{1}{2}\langle\int_0^t\partial_{xx}\varphi(u_s^{\epsilon})-
\varphi''(u^{\epsilon}_s)(\partial_{x}u_s^{\epsilon})^2ds,l\rangle\\
&&=\langle\varphi(u^{\epsilon}_0),l\rangle+\int_0^t\langle\varphi'(u^{\epsilon}_s)l,dW_{s}\rangle
-\frac{1}{2}\langle\int_0^t\varphi''(u_s^{\epsilon})((\partial_{x}u_s^{\epsilon})^2-\int_0^1g_{\epsilon}^2(\cdot,y)dy)ds,
l\rangle\\
&&+\frac{1}{2}\langle\int_0^t\varphi(u_s^{\epsilon}),l''\rangle.
\end{eqnarray*}
Since this identity coincides with the expression derived in
\cite{Z} before taking the limit for $\epsilon\to 0$, the
equivalence of the two formulas is proved.
\item[Gradinaru-Nourdin-Tindel's formula]: In \cite{GNT} the authors
present the following It\^o's formula:
\begin{eqnarray}
\varphi(u_t)=\varphi(u_0)+\int_0^t\langle\varphi'(u_s),\delta
u_s\rangle+\frac{1}{2}\int_0^tTr(e^{2s\Delta}\varphi''(u_s))ds.
\end{eqnarray}
The term $\int_0^t\langle\varphi'(u_s),\delta u_s\rangle$ denotes a
Skorohod's type integral w.r.t. the solution process $u_t(x)$. If we
write formally
$$
\delta u_s(x)=\frac{1}{2}\partial_{xx}u_s(x)ds+\delta W_{s,x}\mbox{
and }\int_0^t\langle\varphi'(u_s),\delta
u_s\rangle=\int_0^t\int_0^1\varphi'(u_s)\diamond
\frac{du_s(x)}{ds}dxds,
$$
we get
\begin{eqnarray*}
\int_0^t\langle\varphi'(u_s),\delta
u_s\rangle&=&\int_0^t\int_0^1\varphi'(u_s)\diamond
\frac{du_s(x)}{ds}dxds\\
&&=\int_0^t\int_0^1\varphi'(u_s)\diamond\Big(\frac{1}{2}\partial_{xx}u_s(x)+
\dot{W}_{s,x}\Big)dxds\\
&&=\int_0^t\int_0^1\varphi'(u_s)\diamond\frac{1}{2}\partial_{xx}u_s(x)dxds+\int_0^t\int_0^1\varphi'(u_s)dW_{s,x}.
\end{eqnarray*}
This procedure, far from being rigorous, suggests some common
feature between (2.1) and (2.4). However the identification of the
"It\^o's terms" of the two formulas, namely
$$
\int_0^t\varphi''(u_s(\cdot))d\sigma^2(s;\cdot)\mbox{ and
}\int_0^tTr(e^{2s\Delta}\varphi''(u_s))ds,
$$
doesn't seem to be straightforward.
\end{description}

\end{document}